\documentclass{amsproc}
\usepackage{amssymb}
\usepackage{graphicx}
\usepackage{amscd}

\theoremstyle{plain}

\newtheorem{corollary}{Corollary}

\newtheorem{remark}{Remark}

\numberwithin{equation}{section}

\input{tcilatex}

\begin{document}
\title[Wavelet frames with Fourier transforms of Laguerre functions]{Wavelet
frames, Bergman spaces and Fourier transforms of Laguerre functions}
\author{Lu\'{\i}s Daniel Abreu}
\urladdr{}
\date{April 14, 2006}
\subjclass{}
\keywords{}
\thanks{}

\begin{abstract}
The Fourier transforms of Laguerre functions play the same canonical role in
wavelet analysis as do the Hermite functions in Gabor analysis. We will use
them as analyzing wavelets in a similar way the Hermite functions were
recently by Gr\"{o}chenig and Lyubarskii in \emph{Gabor frames with Hermite
functions}, C. R. Acad. Sci. Paris, Ser. I 344 157-162 (2007). Building on
Seip%
%TCIMACRO{\U{b4}}%
%BeginExpansion
\'{}%
%EndExpansion
s work \emph{Beurling type density theorems in the unit disc}, Invent.
Math., \textbf{113}, 21-39 (1993), concerning sampling sequences on weighted
Bergman spaces, we find a sufficient density condition for constructing
frames by translations and dilations of the Fourier transform of the $nth$
Laguerre function. As in Gr\"{o}chenig-Lyubarskii theorem, the density
increases with $n$, and in the special case of the hyperbolic lattice in the
upper half plane it is given by%
\begin{equation*}
b\log a<\frac{4\pi }{2n+\alpha },
\end{equation*}%
where $\alpha $ is the parameter of the Laguerre function.
\end{abstract}

\maketitle

\section{Introduction}

The "frame density problem" is one of the fundamental questions in applied
harmonic analysis: how coarse we need to sample a continuous object so that
the resulting discrete object is a frame in a certain Hilbert space? This
question was sharply solved for frames in the Bargmann-Fock \cite{Ly}, \cite%
{S1}, \cite{SW} and in the Bergmann space \cite{Seip2}. However, in the
cases of general wavelet and Gabor transforms, very little is known and only
a few very special windows and analysing wavelets are understood.

The short time Fourier (Gabor) transform with respect to a gaussian window
can be written in terms of the Bargmann transform, mapping isometrically the
space $L^{2}(\mathbf{R})$ onto the Bargmann-Fock space of entire functions.
This is the reason why everything is known about the geometry of sequences
that generate frames by sampling the Gabor transform with gaussian windows $%
g(t)=e^{-\pi t^{2}}$. Apart from this example, the only cases where a
description is known of the lattice sequences that generate frames are the
hyperbolic secant $g(t)=(\cosh at)^{-1}$ \cite{JS} and the characteristic
function of an interval \cite{J2}, which turned out to be a nontrivial
problem. There is also a necessary density condition for Gabor frames due to
Ramanathan and Steger \cite{RS} and it was recently observed that the frame
property is stable under small perturbations of the lattice, given the
window belongs to the Feichtinger algebra \cite{FK}.

If we consider the wavelet case, things appear to be even more mysterious.
Although the density and irregular grid problems have attracted some
attention lately (see \cite{ACM}, \cite{SZ}, \cite{G}, \cite{HK} and
references there in), there is no known counterpart of Ramanathan and
Steeger theorem (there is, however, a version of the HAP property for
wavelet systems \cite{HAP}) and the only information available so far, on
particular systems, concerns a special family of analyzing wavelets that
maps the problem into Bergmann spaces: The wavelet transform, with positive
dilation parameter, with respect to a wavelet of the form (Paul%
%TCIMACRO{\U{b4}}%
%BeginExpansion
\'{}%
%EndExpansion
s Wavelet in some literature) 
\begin{equation}
\psi _{\alpha }(t)=\left( \frac{1}{t+i}\right) ^{\alpha +1},  \label{win}
\end{equation}%
can be rescaled as an isometrical integral transform between spaces of
analytical functions, namely, between the Hardy space on the upper half
plane $H^{2}(\mathbf{U})$ and the Bergman space in the upper half plane $%
A^{2\alpha +1}(\mathbf{U})$. It is natural to refer to this isomorphism as
the Bergman transform (the designation analytic wavelet is also frequently
used). Since the upper half plane can be mapped isomorphically onto the unit
disc by using linear fractional transformations, we can construct a
transform mapping the Hardy space onto the Bergman space in the unit disc.

The approaches used to deal with the special situations mentioned in the
above paragraphs are based in techniques which differ from case to case. It
seems highly desirable to follow a more structured approach. The natural
place to look for this structure is within the context of Hilbert spaces of
analytic functions, where the powerful methods from complex analysis may
answer questions that seem hopeless otherwise.\ Following this line of
reasoning implies using windows that allow to carry the problem to such
spaces. In this direction a major step was taken recently by G\"{o}chenig
and Lyubarskii \cite{GrLy}, by considering Gabor systems with Hermite
functions of order $n$ as windows. They have proved that if the size of the
lattice $\Lambda $ is $<(n+1)^{-1}$ then the referred Gabor system is a
frame and provided an example supporting their conjecture that the result is
sharp.

\ In the Bargmann-Fock setting the Hermite functions play a very special
role \cite[pag. 57]{Gro}. They are, up to normalization constants,
pre-images, under the Bargmann transform, of the monomials $\{z^{n}\}$ and
since the latter constitute an orthogonal basis of the Bargmann-Fock space,
their pre-images also constitute an orthogonal basis of $L^{2}(\mathbf{R})$.
That is, to the isomorphism 
\begin{equation}
L^{2}(\mathbf{R})\overset{\mathcal{B}}{\rightarrow }F^{2}(\mathbf{C})
\label{esquemaGabor1}
\end{equation}
corresponds 
\begin{equation}
h_{n}\overset{\mathcal{B}}{\rightarrow }c_{n}z^{n}  \label{esquemaGabor2}
\end{equation}
where $\mathcal{B}$ is the Bargmann transform, $F^{2}(\mathbf{C})$ is the
Bargmann-Fock space, $h_{n}$ are the Hermite functions and $c_{n}$ some
constants dependent on $n$. They are canonical to time-frequency analysis in
an additional sense, since they constitute the eigenfunctions of the
time-frequency-localization operator with Gaussian window \cite{D1}. The
Hermite functions are eigenfunctions of the Fourier transform and they can
be used \cite{Jan} to describe Feichtinger%
%TCIMACRO{\U{b4}}%
%BeginExpansion
\'{}%
%EndExpansion
s algebra $S_{o}$.

Wavelet and Gabor analysis share many similarities and many of their
structural aspects can be bound together in a more general theory using
representations of locally compact abelian groups \cite{GMP}, \cite{FG}. In
looking for a Wavelet-analogue of Gr\"{o}chenig-Lyubarskii structured
approach to the density problem, we must first clarify what functions should
be used instead of the Hermite functions. In analogy to the above paragraph
it is natural to consider the pre-images, under the Bergman transform, of
the monomials $\{z^{n}\}$ in the unit disc (up to an isomorphism with the
upper half plane). Since the $\{z^{n}\}$ form a basis of the weighted
Bergman spaces in the unit disc (independently of the weight), their
pre-images must be an orthogonal basis of the Hardy space $H^{2}(\mathbf{U})$
and it is reasonable to expect that such functions will play a similar role
in Wavelet analysis as do the Hermite functions in Gabor analysis. Such
functions are the Fourier transforms of the Laguerre functions and they will
be ad-hoc denoted by $S_{n}^{\alpha }$.

Additional structural evidence that the functions $S_{n}^{\alpha }$ are
canonical in Wavelet analysis comes from the work of Daubechies and Paul 
\cite{DP}, where it is shown that they are the eigenfunctions of a
differential operator that commutes with a time--scale localization
operator, once windows of the form (\ref{win}) are chosen. This completely
parallels the situation with the time-frequency-localization operators with
gaussian windows, which commute with the harmonic oscillator and therefore
have as eigenfunctions the Hermite functions \cite{D1}. It was observed by
Seip \cite{Seip0} that these problems have a more natural formulation when
mapped into the convenient spaces of analytical functions.

\section{Description of the results}

There exists an analogue of the correspondence between (\ref{esquemaGabor1})
and (\ref{esquemaGabor2}), but involving four different functional spaces
and the corresponding bases: To the sequence of isomorphisms between the
Hilbert spaces 
\begin{equation}
L^{2}(\mathbf{0,\infty })\overset{\mathcal{F}}{\longleftrightarrow }H^{2}(%
\mathbf{U})\overset{Ber^{\alpha }}{\rightarrow }A_{2\alpha +1}(\mathbf{U})%
\overset{T_{\alpha }}{\rightarrow }A_{2\alpha +1}(\mathbf{D})\text{,}
\label{esquemawavelet1}
\end{equation}
where $A_{2\alpha +1}(\mathbf{U})$ and $A_{2\alpha +1}(\mathbf{D})$\ are the
weighted Bergmann spaces in the upper-half plane and in the unit disc,
respectively, corresponds the relations between the basis of the respective
spaces:

\begin{equation}
l_{n}^{\alpha }\overset{\mathcal{F}}{\longleftrightarrow }S_{n}^{\alpha }%
\overset{Ber^{\alpha }}{\rightarrow }c_{n}^{\alpha }\Psi _{n}^{\alpha }%
\overset{T_{\alpha }}{\rightarrow }c_{n}^{\alpha }z^{n}\text{.}
\label{esquemawavelet2}
\end{equation}
for some constants $c_{n}^{\alpha }$, where $l_{n}^{\alpha }$ is the $nth$
Laguerre function of order $n$ and parameter $\alpha $, $S_{n}^{\alpha }$ is
its Fourier transform and $\Psi _{n}^{\alpha }$ is a basis of $A_{2\alpha
+1}(\mathbf{U})$ to be defined in section 4. This correspondence was
implicit in the conection between papers \cite{DP} and \cite{Seip0} but
since it was not stated explicitly, we devote section 4 to clarify how
exactly does it work.

The functions $S_{n}^{\alpha }$ were computed recently in closed form \cite%
{Shen}, but the connection to the Wavelet transform seems to have been
unnoticed. It is also shown in \cite{Shen} that $S_{n}^{\alpha }$ are, up to
a fractional transformation, defined in terms of a certain system of
orthogonal polynomials on the unit circle. We will show that this family of
polynomials is nothing more but the circular Jacobi orthogonal polynomials
for which there is, for computational purposes, a very convenient three term
recurrence formula. This connection will make our case to call the functions 
$S_{n}^{\alpha }$ the \emph{rational Jacobi orthogonal functions}.

With a computational method available for evaluating the functions $%
S_{n}^{\alpha }$ and graphic evidence of their good localization properties
(see for example the plots of their real versions in page 44 of \cite{Dau}),
it is natural to investigate how to use such functions as analysing wavelets
and to obtain frames from the resulting discretization. We will obtain the
following sufficient condition on the density of the parameters of the
hyperbolic lattice $\{(a^{j}bk,a^{j})\}_{j,k\in \mathbf{Z}}$: 
\begin{equation*}
b\log a<\frac{4\pi }{2n+\alpha }
\end{equation*}%
that is, as in Gr\"{o}chenig and Lyubarskii result \cite{GrLy}, this density
increases with $n$. This will follow as a special case of the following more
general theorem, which uses Seip%
%TCIMACRO{\U{b4}}%
%BeginExpansion
\'{}%
%EndExpansion
s notion of lower Beurling density in the unit disk and constitutes the main
result in this paper:

\textbf{Theorem 5.1 }\emph{Let }$\Gamma \subset D$\emph{\ denote a separated
sequence obtained from mapping the sequence }$\{z_{m,j}=a_{m,j}+ib_{m,j}\}%
\subset U$\emph{\ into the unit disk via a Cayley transform. If }$%
D^{-}(\Gamma )>n+\frac{\alpha }{2}$\emph{\ then }$%
\{T_{a_{m,j}}D_{b_{m,j}}S_{n}^{\alpha }(\frac{t}{2})\}_{m,j}$\emph{\ is a
frame of }$H^{2}(U).$

Case $\alpha =0$ gives the sufficient part of theorem 7.1 in \cite{Seip2}.\
The proof of theorem 1 will be dramatically simplified with the observation
of the following remarkable ocurrence: There exist constants $\{C\}$ and
sequences of numbers $\{a_{k}\}$, both depending on $\alpha $ and $n$, such
that 
\begin{equation}
S_{n}^{\alpha }(\frac{t}{2})=C\sum_{k=0}^{n}a_{k}\psi _{k+\frac{\alpha }{2}%
}(t)\text{.}  \label{repcomb}
\end{equation}%
This will allow to obtain a formula which expresses the wavelet transform
with window $S_{n}^{\alpha }$ as a linear combination of simple but
non-analytic functions. As a result we will have to deal with a situation
that is reminiscent of the non-analiticity problem that Gr\"{o}chenig and
Lyubarskii faced with formula (15) of \cite{GrLy} and which they solved by
using the so called Wexler-Raz identities \cite{Gro}. Since there is no
known analogue of the Wexler-Raz identities in wavelet analysis, we have to
deal with this in a more direct way. The key idea is to map the problem into
the unit circle and use a deep sampling formula proved by Seip in \cite%
{Seip2}. Althought this approach lacks the simplicity of \cite{GrLy} it
gives the result for more general sequences, while the method in \cite{GrLy}
works only for the lattice case due to the restritions imposed by the use of
the Wexler-Raz identities.

We organize our ideas as follows. The next section contains the main
definitions and facts concerning wavelet transforms, Bergman spaces and
Laguerre functions. Section four contains the evaluation of the pre-images
required to build up the correspondences (\ref{esquemawavelet1}) and (\ref%
{esquemawavelet2}). The fifth section contains our main results on wavelet
frames with Fourier transforms of Laguerre functions. We conclude the paper
collecting some further properties of the functions $S_{n}^{\alpha }$ and
clarifying their classification within known families of orthogonal
polynomials.

\section{Tools}

\subsection{The Bergman transform}

Now we present a sinthesis of ideas that appeared in the section 3.2 of \cite%
{GMP} and in \cite{DP} (see also \cite[pag. 31]{Dau}) Here they will be
exposed in such a way that the role of the Bergman spaces is emphasized.

\ Consider the dilation and translation operators 
\begin{eqnarray*}
D_{s}f(x) &=&\left| s\right| ^{-\frac{1}{2}}f(s^{-1}x) \\
T_{x}f(t) &=&f(t-x)
\end{eqnarray*}
and define 
\begin{equation*}
\psi ^{x,s}(t)=T_{x}D_{s}\psi (t)=\left| s\right| ^{-\frac{1}{2}}\psi (\frac{%
t-x}{s})\text{.}
\end{equation*}
The wavelet transform of a function $f$, with respect to the wavelet $\psi $
is 
\begin{equation*}
W_{\psi }f(x,s)=\left\langle f,T_{x}D_{s}\psi (t)\right\rangle _{L^{2}(%
\mathbf{R})}=\int_{-\infty }^{\infty }f(t)\overline{\psi ^{x,s}(t)}dt\text{.}
\end{equation*}
A function $\psi \in L^{2}(R)$ is said to be admissible if 
\begin{equation*}
\int_{0}^{\infty }\left| \mathcal{F}\psi (s)\right| ^{2}\frac{ds}{s}=K
\end{equation*}
where $K$ is a constant. If $\psi $ is admissible, then for all $f\in
L^{2}(R)$ we have 
\begin{equation}
\int_{0}^{\infty }\int_{-\infty }^{\infty }s^{-2}\left| W_{\psi
}f(x,s)\right| ^{2}dxds=K\left\| f\right\| ^{2}  \label{isometry}
\end{equation}

We will restrict ourselves to parameters $s>0$ and functions $f\in H^{2}(%
\mathbf{U})$, where $H^{2}(\mathbf{U})$ is the Hardy space in the upper half
plane $\mathbf{U}=\{z=x+is:s>0\}$%
\begin{equation*}
H^{2}(\mathbf{U})=\{f:f\text{ is analytic in }\mathbf{U}\text{ and }%
\sup_{0<s<\infty }\int_{-\infty }^{\infty }\left| f(x+is)\right|
^{2}dx<\infty \}.
\end{equation*}
Let $\mathcal{F}$ denote the Fourier transform 
\begin{equation*}
(\mathcal{F}f)(t)=\frac{1}{\sqrt{2\pi }}\int_{-\infty }^{\infty
}e^{-itx}f(x)dx
\end{equation*}
By the Paley-Wiener theorem,\ $H^{2}(\mathbf{U})$\ is constituted by the
functions whose Fourier transform is supported in $\left( 0,+\infty \right) $
and belongs to $L^{2}(0,\infty )$.

Now take the special window $\psi _{\alpha }(t)$ defined in (\ref{win}).
Since 
\begin{equation}
\mathcal{F}\overline{\psi _{\alpha }^{-x,s}(t)}=\mathbf{1}_{\left[ 0,\infty %
\right] }s^{\alpha +\frac{1}{2}}t^{\alpha }e^{i(x+is)t}  \label{Fg}
\end{equation}%
then 
\begin{equation}
W_{\psi _{\alpha }}f(-x,s)=s^{\alpha +\frac{1}{2}}\int_{0}^{\infty
}t^{\alpha }e^{izt}(\mathcal{F}f)(t)dt  \label{anwav}
\end{equation}%
where the function defined by the integral is analytic in $z=x+si$. The
identity (\ref{isometry}) gives 
\begin{equation}
\int \int_{\mathbf{U}}\left\vert W_{\psi _{\alpha }}f(-x,s)\right\vert
^{2}s^{-2}dxds=\left\Vert f\right\Vert _{H^{2}(\mathbf{R})}  \label{estwav}
\end{equation}%
This motivates the definition of the \emph{Bergman transform}, or the \emph{%
analytic wavelet transform}: 
\begin{equation}
Ber^{\alpha }\text{ }f(z)=\int_{0}^{\infty }t^{\alpha }e^{izt}(\mathcal{F}%
f)(t)dt=s^{-\alpha -\frac{1}{2}}W_{\psi _{\alpha }}f(-x,s),  \label{Ber}
\end{equation}%
where $z=x+is$ (see for instance \cite{HK}, where the authors use this
Bergman transform with the same normalization in the case $\alpha =1$).
Introducing the scale of weighted Bergman spaces 
\begin{equation*}
A_{\alpha }(\mathbf{U})=\{f\text{ analytic in }\mathbf{U}\text{ such that }%
\int \int_{\mathbf{U}}\left\vert f(z)\right\vert ^{2}s^{\alpha
-2}dxds<\infty \}\text{,}
\end{equation*}%
it is clear from (\ref{estwav}) and (\ref{Ber}) that\ $Ber^{\alpha }$ $%
f(z)\in A_{2\alpha +1}(\mathbf{U})$. We have therefore an isometric
transformation 
\begin{equation*}
Ber^{\alpha }:H^{2}(\mathbf{R})\rightarrow A_{2\alpha +1}(\mathbf{U})
\end{equation*}%
The weighted Bergman spaces in the unit disc are denoted by $A_{\alpha }(%
\mathbf{D})$ and defined as 
\begin{equation*}
A_{\alpha }(\mathbf{D})=\left\{ f\text{ analytic in }D\text{ such that }\int
\int_{D}\left\vert f(z)\right\vert ^{2}(1-\left\vert z\right\vert )^{\alpha
-2}dxdy<\infty \right\} \text{.}
\end{equation*}%
For detailed expositions of the theory of Bergman spaces we point out the
monographs \cite{Z}, \cite{DS}, \cite{HKZ}

\subsection{Frames and sampling sequences}

A sequence of functions $\{e_{j}\}$ is said to be a frame in a Hilbert space 
$H$ if there exist constants $A$ and $B$ such that 
\begin{equation}
A\left\Vert f\right\Vert ^{2}\leq \sum_{j}\left\vert \left\langle
f,e_{j}\right\rangle \right\vert ^{2}\leq B\left\Vert f\right\Vert ^{2}.
\label{frame}
\end{equation}

In spaces of analytic functions a related concept to frames is the one of a
sampling sequence. A set $\Gamma =\{z_{j}\}$ is said to be a sampling
sequence for the Bergman space $A_{\alpha }(\mathbf{U})$ if there exist
positive constants $A$ and $B$ such that 
\begin{equation}
A\int \int_{\mathbf{U}}\left\vert f(z)\right\vert ^{2}y^{\alpha -2}dxds\leq
\sum_{j}\left\vert f(z_{j})\right\vert ^{2}y_{j}^{\alpha }\leq B\int \int_{%
\mathbf{U}}\left\vert f(z)\right\vert ^{2}y^{\alpha -2}dxds.
\label{sampling}
\end{equation}

The corresponding definition of a sampling sequence for the Bergman space on
the unit disc is%
\begin{eqnarray*}
A\int \int_{\mathbf{D}}\left\vert f(z)\right\vert ^{2}(1-\left\vert
z\right\vert ^{2})^{\alpha -2}dz &\leq &\sum_{j}\left\vert
f(z_{j})\right\vert ^{2}(1-\left\vert z\right\vert ^{2})^{\alpha } \\
&\leq &B\int \int_{\mathbf{D}}\left\vert f(z)\right\vert ^{2}(1-\left\vert
z\right\vert ^{2})^{\alpha -2}dz.
\end{eqnarray*}%
Define the pseudohyperbolic metric on the unit disk by 
\begin{equation*}
\varrho (z,\zeta )=\left\vert \frac{z-\zeta }{1-\overline{\zeta }z}%
\right\vert
\end{equation*}%
Following \cite{Seip2}, a sequence $\Gamma =\{z_{j}\}\subset \mathbf{D}$ is
separated if $\inf_{n\neq j}\left\vert \frac{z_{j}-z_{n}}{z_{j}-\overline{z}%
_{n}}\right\vert >0$ and its lower density $D^{-}$ is given by 
\begin{equation*}
D^{-}(\Gamma )=\lim_{r\rightarrow 1}\inf \inf_{z}\frac{\sum_{\varrho
(z_{j},z)<r}(1-\varrho (z_{j},z))}{\log \frac{1}{1-r}}
\end{equation*}%
The next Theorem is from \cite{Seip2}:

\textbf{Theorem A }\emph{A separated sequence }$\Gamma \subset \mathbf{D}\ 
\emph{is}$ $\emph{a}$ $\emph{sampling}$ $\emph{sequence}$ $\emph{for}$ $%
A_{\alpha }(\mathbf{D})$ iff $D^{-}(\Gamma )>\frac{\alpha }{2}-\frac{1}{2}$.

It is a generalization of an earlier lattice result \cite{Seip}:

\textbf{Theorem B }\emph{Let }$\Gamma (a,b)=\{z_{jk}\}_{j,k\in \mathbf{Z}}$%
\emph{, where}$\ z_{jk}=a^{j}(bk+i)$\emph{. }$\Gamma (a,b)$\emph{\ is a
sampling sequence for }$A_{\alpha }(\mathbf{U})$\emph{\ if and only if }$%
b\log a<\frac{4\pi }{\alpha -1}$\emph{.}

\subsection{Fourier transforms of Laguerre functions}

The Laguerre polynomials will play a central role in our discussion. One way
to define them is by means of the Rodrigues formula 
\begin{equation}
L_{n}^{\alpha }(x)=\frac{e^{x}x^{-\alpha }}{n!}\frac{d^{n}}{dx^{n}}\left[
e^{-x}x^{\alpha +n}\right]  \label{Rodr}
\end{equation}%
and this gives, in power series form 
\begin{equation*}
L_{n}^{\alpha }(x)=\frac{(\alpha +1)_{n}}{n!}\sum_{k=0}^{n}\frac{(-n)_{k}}{%
(\alpha +1)_{k}}\frac{x^{k}}{k!}\text{.}
\end{equation*}%
For information on specific systems of orthogonal polynomials, we suggest 
\cite{Ism}.\ The Laguerre functions are defined as 
\begin{equation*}
l_{n}^{\alpha }(x)=\mathbf{1}_{\left[ 0,\infty \right] }(x)e^{-x/2}x^{\alpha
/2}L_{n}^{\alpha }(x)
\end{equation*}%
and they are known to constitute an orthogonal basis for the space $%
L^{2}(0,\infty )$. By the Paley-Wiener theorem, the Fourier transform is an
isomorphism between $H^{2}(\mathbf{U})$\ and $L^{2}(0,\infty )$. Therefore,
the Fourier transform of the functions $l_{n}^{\alpha }$ form an orthogonal
basis for the space $H^{2}(\mathbf{U})$. From now on we will set 
\begin{equation}
\mathcal{F}S_{n}^{\alpha }(t)=l_{n}^{\alpha }(t)  \label{Sn}
\end{equation}%
The functions $S_{n}^{\alpha }(t)$ can be evaluated explictly using (\ref%
{Rodr}). This was already done by Shen in \cite{Shen}, where the author was
interested in describing the orthogonal polynomials which arise from an
application of the Fourier transform on the Laguerre polynomials. A
description of this method of generating new families of orthogonal
polynomials with an interesting historical account is in the paper \cite{Koe}
where a link is provided between Jacobi and Meixner-Polaczek polynomials.
Following \cite{Shen}, we have

\begin{equation}
S_{n}^{\alpha }(t)=\frac{\Gamma \left( \frac{\alpha }{2}+1\right) (1+\alpha
)_{n}}{n!}\sum_{k=0}^{n}\frac{(-n)_{k}(\frac{\alpha }{2}+1)_{k}}{k!(\alpha
+1)_{k}}\left( \frac{1}{\frac{1}{2}-it}\right) ^{k+\frac{\alpha }{2}+1}
\label{Series}
\end{equation}

\section{Bergman transform of $S_{n}^{\protect\alpha }(t)$}

In this section the correspondences (\ref{esquemawavelet1}) and (\ref%
{esquemawavelet2}) will be established via a direct calculation. The complex
linear fractional transformations will play an important role, in a style
that is reminiscent of the way they are used in the discrete series
representation of $SL(2,\mathbf{R})$ over the Bergman space \cite[chapter IX]%
{Lang}.

We first define a set of functions that, as we shall see later, constitute a
basis of $A_{\alpha }(\mathbf{U}).$ For every $n\geq 0$ and $\alpha >-1$ let 
\begin{equation*}
\Psi _{n}^{\alpha }(z)=\left( \frac{iz+\frac{1}{2}}{iz-\frac{1}{2}}\right)
^{n}\left( iz-\frac{1}{2}\right) ^{-\alpha -1}\text{.}
\end{equation*}

Now define a map $T_{\alpha }$ such that for every function $f\in A_{\alpha
}($\textbf{$U$}$)$ the action of $T_{\alpha }$ is 
\begin{equation}
T_{\alpha }f(w)=f\left[ \frac{i}{2}\frac{w+1}{w-1}\right] \left( \frac{1}{1-w%
}\right) ^{\alpha +1}  \label{Talpha}
\end{equation}

The range space of $T_{\alpha }$\ is a weighted Bergman space in the unit
disc.

\textbf{Lemma 4.1 }\emph{The map } 
\begin{equation*}
T_{\alpha }:A_{\alpha }(\mathbf{U})\rightarrow A_{\alpha }(\mathbf{D})
\end{equation*}
\emph{is an unitary isometry between Hilbert spaces.}

\textbf{Proof. }Argue as in the proof of lemma 1 in \cite[pag.185]{Lang}. $%
\Box $

\textbf{Proposition 4.1 }\emph{For }$n=0,1,...$\emph{, the following
relations hold: } 
\begin{equation}
Ber^{\alpha }(S_{n}^{2\alpha })=c_{n}^{\alpha }\Psi _{n}^{\alpha }\text{ },%
\text{ with\ }c_{n}^{\alpha }=(-1)^{\alpha +1}(2\alpha +n)!/n!  \label{tr1}
\end{equation}%
\emph{and, for }$\left\vert z\right\vert <1$\emph{, } 
\begin{equation}
T_{\alpha }(\Psi _{n}^{\alpha })=z^{n}  \label{tr2}
\end{equation}%
\emph{In other words, the function }$S_{n}^{2\alpha }$\emph{\ is the
pre-image, under }$T_{\alpha }\circ \left( \frac{1}{c_{n}^{\alpha }}%
Ber^{\alpha }\right) $\emph{, of }$z^{n}\in A_{\alpha }(D)$\emph{.}

\textbf{Proof. }Using the definition of the Bergman transform (\ref{Ber})
and (\ref{Sn}) we have 
\begin{eqnarray*}
Ber^{\alpha }\text{ }S_{n}^{2\alpha }(z) &=&\int_{0}^{\infty }t^{\alpha
}e^{izt}(\mathcal{F}S_{n}^{2\alpha }(t))(t)dt \\
&=&\int_{0}^{\infty }e^{(iz-\frac{1}{2})t}t^{2\alpha }L_{n}^{2\alpha }(t)dt
\\
&=&\frac{1}{n!}\mathcal{L}\left[ \frac{d^{n}}{dx^{n}}\left[ e^{-x}x^{2\alpha
+n}\right] \right] (-iz+\frac{1}{2}) \\
&=&c_{n}^{\alpha }\Psi _{n}^{2\alpha }(z)\text{,}
\end{eqnarray*}%
where in the last two identities we are applying Rodrigues formula (\ref%
{Rodr}) for the Laguerre polynomials and writing the integral in terms of
the Laplace transform $\mathcal{L}$, whose well known properties establish
the last identity. Now, the linear fractional transformation 
\begin{equation*}
w=\frac{iz+\frac{1}{2}}{iz-\frac{1}{2}}
\end{equation*}%
is an analytic isomorphism between the upper half plane and the unit circle.
Since the inverse of this transformation is given by 
\begin{equation*}
z=\frac{i}{2}\frac{w+1}{w-1}\text{,}
\end{equation*}%
a short calculation with the definition of $T_{\alpha }$ gives(\ref{tr2}). $%
\Box $

\begin{remark}
$\{\Psi _{n}^{\alpha }(z)\}$ is a basis of $A_{\alpha }($\textbf{$U$}$)$ and
the map $T_{\alpha }$ is an unitary isomorphism between $A_{\alpha }($%
\textbf{$U$}$)$ and $A_{\alpha }(\mathbf{D})$. Indeed, since $\{z^{n}\}$ is
an orthogonal basis of the space $A_{\alpha }(\mathbf{D})$ \cite[pag. 186]%
{Lang} and it is contained in the range of $T_{\alpha }$, $T_{\alpha }$ is
onto and therefore an unitary isomorphism.The functions $\{\Psi _{n}^{\alpha
}(z)\}$ form a basis for the space $A_{\alpha }(\mathbf{U})$ since they are
the pre-images of the basis $\{z^{n}\}$.
\end{remark}

As a consequence we obtain a new proof of the (known) isomorphic property of
the Bergman transform.

\begin{corollary}
\ The transform $Ber^{\alpha }:H^{2}(\mathbf{R})\rightarrow A_{2\alpha +1}(%
\mathbf{U})$\ is an isometric isomorphism.
\end{corollary}

\textbf{Proof. }The isometry is a consequence of the isometric property of
the wavelet transform, so we need only to prove that the $Ber^{\frac{\alpha 
}{2}}$ is onto. But in view of the preceding remark and Theorem 5.1, the
range of $Ber^{\frac{\alpha }{2}}$ contains a basis of $A_{\alpha }(\mathbf{U%
})$. Therefore, $Ber^{\frac{\alpha }{2}}$ is onto. $\Box $

\begin{remark}
Similar calculations as we have seen here also play a role in \cite{DOZ}, in
the context of Laplace transformations and group representations and in \cite%
{DP}, to obtain an explicit formula for the eigenvalues of the time-scale
localization operator.
\end{remark}

\section{Wavelet frames with Fourier transforms of Laguerre functions}

\subsection{Preparation}

We wish to construct wavelet frames with analysing wavelets $S_{n}^{\alpha }$%
, by using the common discretization of the continuous wavelet transform via
the hyperbolic lattice. Discretizing the scale parameter $s$ in the Wavelet
transform by a sequence $a^{j}$ and the parameter $x$ by $a^{j}bm$ gives 
\begin{equation}
W_{\psi }f(a^{j}bm,a^{j})=\left\langle f,T_{a^{j}bm}D_{a^{j}}\psi
\right\rangle .  \label{sampl}
\end{equation}%
We want to know conditions in $a$ and $b$ under which $%
\{T_{a^{j}bk}D_{a^{j}}\psi \}$ is a frame, for a given wavelet $\psi $. More
generally we can consider the non lattice problem replacing $(a^{j}bm,a^{j})$
by more sequences $(a_{jm},b_{jm})$.

\begin{remark}
Observe that it follows immeadiately from Theorem A and Theorem B that if $%
\Gamma \subset \mathbf{D}$ is the sequence obtained from mapping the
sequence $\{z_{k,j}=a^{j}(bk+i)\}\subset \mathbf{U}$ into the unit disk via
a Cayley transform, then $D^{-}(\Gamma )=\frac{2\pi }{b\log a}$ (to move
between the notation of \cite{Seip2} and \cite{Seip} set $\alpha =2n+1$).
\end{remark}

\begin{remark}
Observe that $\{T_{a^{j}bm}D_{a^{j}}\psi _{\alpha }\}_{j,k\in \mathbf{Z}}$
is a frame of $H^{2}(\mathbf{U})$ iff $\Gamma
(a,b)=\{z=a^{j}bm+a^{j}i\}_{j,k\in \mathbf{Z}}$ is a set of sampling for $%
A_{2\alpha +1}(\mathbf{U}).$ Indeed, since functions in $A_{2\alpha +1}(%
\mathbf{U})$ can be identified with $Ber^{\alpha }$ transforms of $H^{2}(%
\mathbf{U})$ functions, it follows from (\ref{sampling}) and (\ref{Be}) that 
$\Gamma (a,b)$ is a set of sampling for $A_{2\alpha +1}(\mathbf{U})$ if and
only if 
\begin{equation*}
A\int \int_{\mathbf{U}}\left\vert W_{\psi _{\alpha }}f(x,s)\right\vert ^{2}%
\frac{dxds}{s^{2}}\leq \sum_{j,m}\left\vert W_{\psi _{\alpha
}}f(a^{j}bm,a^{j})\right\vert ^{2}\leq B\int \int_{\mathbf{U}}\left\vert
W_{\psi _{\alpha }}f(x,s)\right\vert ^{2}\frac{dxds}{s^{2}}\text{.}
\end{equation*}%
Using (\ref{estwav}) this is equivalent to 
\begin{equation*}
A\left\Vert f\right\Vert _{H^{2}(\mathbf{U})}^{2}\leq \sum_{j,m}\left\vert
\left\langle f,T_{a^{j}bm}D_{a^{j}}\psi _{\alpha }\right\rangle \right\vert
^{2}\leq B\left\Vert f\right\Vert _{H^{2}(\mathbf{U})}^{2}
\end{equation*}%
which says that $\{T_{a^{j}bm}D_{a^{j}}\psi _{\alpha }\}_{j,k}$ is a frame
in $H^{2}(\mathbf{U})$.
\end{remark}

We now collect the required preliminary lemmas for the proof of the main
result. The first was quoted in the introduction and is just a simple
modification of the representation (\ref{Series}). The second appeared in 
\cite{Seip} (see \cite[pag. 160]{DS} for a proof) and the third one is a
quite deep sampling theorem which is formula (30) in \cite{Seip2} with $s=0$
and $\epsilon =0$ (see also \cite[pag. 216]{DS} for a direct derivation with
this choice of parameters).

\textbf{Lemma 5.2 }\emph{The functions }$S_{n}^{\alpha }$\emph{\ can written
as the linear combination (\ref{repcomb}) of analyzing wavelets }$\psi _{k+%
\frac{\alpha }{2}}(t)$\emph{\ defined by }%
\begin{equation*}
S_{n}^{\alpha }(\frac{t}{2})=C\sum_{k=0}^{n}a_{k}\psi _{k+\frac{\alpha }{2}%
}(t)\text{,}
\end{equation*}%
\emph{\ where the constants }$C$\emph{\ and the coeficients }$a_{k}$\emph{\
given as } 
\begin{equation*}
C=C^{\alpha ,n}=\frac{\Gamma \left( \frac{\alpha }{2}+1\right) (1+\alpha
)_{n}}{n!}
\end{equation*}%
\begin{equation*}
a_{k}=a_{k}^{\alpha ,n}=(2i)^{k+\frac{\alpha }{2}+1}\frac{(-n)_{k}(\frac{%
\alpha }{2}+1)_{k}}{k!(\alpha +1)_{k}}.
\end{equation*}

\textbf{Lemma 5.3 } \emph{For }$1<s<t$ \emph{we have}%
\begin{equation}
\sum_{j}\frac{(1-\left\vert z_{j}\right\vert ^{2})^{s}}{\left\vert 1-%
\overline{\zeta }z_{j}\right\vert ^{t}}\leq C(1-\left\vert \zeta \right\vert
^{2})^{s-t}  \label{soma}
\end{equation}

\emph{for any separated sequence }$\{z_{j}\}$.

\textbf{Lemma 5.4 }\emph{Let }$D(\{w_{j}\})>\beta $. \emph{Then every }$f\in
A_{2\beta +1}(D)$\emph{\ satisfies}%
\begin{equation}
(1-\left\vert w\right\vert ^{2})^{\beta }f(w)=\sum_{j}(1-\left\vert
w_{j}\right\vert ^{2})^{\beta }f(w_{j})h_{j}(w)\text{,}
\label{samplingformula}
\end{equation}%
\emph{with the estimate for }$h_{j}(w)$\emph{\ }%
\begin{equation}
\left\vert h_{j}(w)\right\vert \leq C\frac{(1-\left\vert w\right\vert
^{2})(1-\left\vert w_{j}\right\vert ^{2})}{\left\vert 1-\overline{w_{j}}%
w\right\vert ^{2}}\text{.}  \label{h2}
\end{equation}

\subsection{Main result}

To prove our main result we will adhere to the following plan: we begin by
deriving from Lemma 5.2 a formula expressing the wavelet transform with
analysing wavelet $S_{n}^{\alpha }$ in terms of the wavelet transform with
analysing wavelet $\psi _{k+\frac{\alpha }{2}}$. The proof of the right
frame inequality is then straightforward (as usual in frame inequalities).
To prove the left hand inequality we first rewrite the formula in terms of
complex functions defined on the upper half plane and then map it into the
unit disk by means of a Cayley transform. Manipulating the resulting formula
and using lemma 5.3 and lemma 5.4 gives an estimate in the unit disk. We
proceed backwards to the half plane to end up with the required estimate in
terms of wavelet frames.%
\begin{equation*}
\end{equation*}

\textbf{Theorem 5.1 }\emph{Let }$\Gamma \subset D$\emph{\ denote a separated
sequence obtained from mapping the sequence }$\{z_{m,j}=a_{m,j}+ib_{m,j}\}%
\subset U$\emph{\ into the unit disk via a Cayley transform. If }$%
D^{-}(\Gamma )>n+\frac{\alpha }{2}$\emph{\ then }$%
\{T_{a_{m,j}}D_{b_{m,j}}S_{n}^{\alpha }(\frac{t}{2})\}_{m,j}$\emph{\ is a
frame of }$H^{2}(U).$

\textbf{Proof. }The definition of the wavelet transform gives, using (\ref%
{repcomb}), 
\begin{eqnarray*}
W_{S_{n}^{\alpha }(\frac{t}{2})}f(-x,s) &=&\left\langle
f,T_{-x}D_{s}S_{n}^{\alpha }(\frac{t}{2})\right\rangle \\
&=&C\sum_{k=0}^{n}a_{k}\left\langle f,T_{-x}D_{s}\psi _{k+\frac{\alpha }{2}%
}\right\rangle \\
&=&C\sum_{k=0}^{n}a_{k}W_{\psi _{k+\frac{\alpha }{2}}}f(-x,s)\text{.}
\end{eqnarray*}

The right frame inequality is now easy: observe that, since $%
\{T_{a_{m,j}}D_{b_{m,j}}\psi _{k+\frac{\alpha }{2}}\}$ is a frame, then for
every $k=1,...,n$ there exists a $B_{k}$ such that%
\begin{equation*}
\sum_{j,m}\left\vert \left\langle f,T_{a_{m,j}}D_{b_{m,j}}\psi _{k+\frac{%
\alpha }{2}}\right\rangle \right\vert ^{2}\leq B_{k}\left\Vert f\right\Vert
_{H^{2}(\mathbf{U})}^{2}.
\end{equation*}%
Let $B_{n}=\max_{1\leq k\leq n}B_{k}$. Then, 
\begin{eqnarray*}
\sum_{j,m}\left\vert \left\langle f,T_{a_{m,j}}D_{b_{m,j}}S_{n}^{\alpha }(%
\frac{t}{2})\right\rangle \right\vert ^{2} &=&\sum_{j,m}\left\vert
C\sum_{k=0}^{n}a_{k}\left\langle f,T_{a_{m,j}}D_{b_{m,j}}\psi _{k+\frac{%
\alpha }{2}}\right\rangle \right\vert ^{2} \\
&\leq &C\sum_{k=0}^{n}\left\vert a_{k}\right\vert ^{2}\sum_{j,m}\left\vert
\left\langle f,T_{a_{m,j}}D_{b_{m,j}}\psi _{k+\frac{\alpha }{2}%
}\right\rangle \right\vert ^{2} \\
&\leq &B\left\Vert f\right\Vert ^{2}\text{,}
\end{eqnarray*}%
where $B=B_{n}C\sum_{k=0}^{n}\left\vert a_{k}\right\vert ^{2}$ and we are
done.

Now the left frame inequality. Define the function $F(z)$ as 
\begin{equation}
F(z)=W_{S_{n}^{\alpha }(\frac{t}{2})}f(-x,s)=C\sum_{k=0}^{n}a_{k}s^{k+\frac{%
\alpha }{2}+\frac{1}{2}}Ber^{k+\frac{\alpha }{2}}f(z)  \label{Formula}
\end{equation}%
where $z=x+is$ belongs to the upper half plane $\mathbf{U}$.\ We construct a
related function in the unit disc, by setting $w=\frac{z-i}{z+i}$:%
\begin{equation*}
G(w)=C\sum_{k=0}^{n}a_{k}(1-\left\vert w\right\vert ^{2})^{k+\frac{\alpha }{2%
}+\frac{1}{2}}G_{k+\frac{\alpha }{2}}(w)
\end{equation*}%
where 
\begin{equation*}
G(w)=F(i\frac{1+w}{1-w})
\end{equation*}%
and%
\begin{equation*}
G_{k+\frac{\alpha }{2}}(w)=\left\vert \frac{2}{1-w}\right\vert ^{2k+\alpha
+1}Ber^{k+\frac{\alpha }{2}}f(i\frac{1+w}{1-w})\in A_{2k+\alpha +1}(\mathbf{D%
}).
\end{equation*}%
\ Therefore, if $D(\{w_{j}\})>n+\frac{\alpha }{2}$ we can use (\ref%
{samplingformula}) to write, for every $k\leq n$,%
\begin{equation*}
(1-\left\vert w\right\vert ^{2})^{k+\frac{\alpha }{2}}G_{k+\frac{\alpha }{2}%
}(w)=\sum_{j}(1-\left\vert w_{j}\right\vert ^{2})^{k+\frac{\alpha }{2}}G_{k+%
\frac{\alpha }{2}}(w_{j})h_{j}(w)\text{.}
\end{equation*}%
As a result,%
\begin{eqnarray*}
(1-\left\vert w\right\vert ^{2})^{-\frac{1}{2}}G(w) &=&C^{\alpha
,n}\sum_{k=0}^{n}a_{k}(1-\left\vert w\right\vert ^{2})^{k+\frac{\alpha }{2}%
}G_{k+\frac{\alpha }{2}}(w) \\
&=&C\sum_{k=0}^{n}a_{k}\sum_{j}(1-\left\vert w_{j}\right\vert ^{2})^{k+\frac{%
\alpha }{2}}G_{k+\frac{\alpha }{2}}(w_{j})h_{j}(w) \\
&=&\sum_{j}h_{j}(w)C\sum_{k=0}^{n}a_{k}(1-\left\vert w_{j}\right\vert
^{2})^{k+\frac{\alpha }{2}}G_{k+\frac{\alpha }{2}}(w_{j}) \\
&=&\sum_{j}h_{j}(w)(1-\left\vert w_{j}\right\vert ^{2})^{-\frac{1}{2}%
}G(w_{j}).
\end{eqnarray*}%
Now use estimate (\ref{h2}), then Cauchy-Schwarz and finally (\ref{soma}) in
the following way.%
\begin{eqnarray*}
(1-\left\vert w\right\vert ^{2})^{-1}\left\vert G(w)\right\vert ^{2} &\leq & 
\left[ C\sum_{j}\frac{(1-\left\vert w\right\vert ^{2})}{\left\vert 1-%
\overline{w_{j}}w\right\vert ^{2}}(1-\left\vert w_{j}\right\vert ^{2})^{%
\frac{1}{2}}G(w_{j})\right] ^{2} \\
&\leq &C\sum_{j}(1-\left\vert w\right\vert ^{2})^{2}\left\vert
G(w_{j})\right\vert ^{2}\sum_{j}\frac{(1-\left\vert w_{j}\right\vert
^{2})^{1}}{\left\vert 1-\overline{w_{j}}w\right\vert ^{4}} \\
&\leq &C\sum_{j}(1-\left\vert w\right\vert ^{2})^{-1}\left\vert
G(w_{j})\right\vert ^{2}.
\end{eqnarray*}%
Thus the factor $(1-\left\vert w\right\vert ^{2})^{-1}$ cancels and we have%
\begin{equation*}
\left\vert G(w)\right\vert ^{2}\leq C\sum_{j}\left\vert G(w_{j})\right\vert
^{2}\text{.}
\end{equation*}%
This last inequality gives at once%
\begin{equation*}
\int \int_{D}\left\vert G(w)\right\vert ^{2}(1-\left\vert w\right\vert
^{2})^{-2}dw\leq C\sum_{j}\left\vert G(w_{j})\right\vert ^{2}
\end{equation*}%
or, in the upper half plane,%
\begin{equation*}
A\int \int_{\mathbf{U}}\left\vert F(z)\right\vert ^{2}s^{-2}dxds\leq
\sum_{j}\left\vert F(z_{j})\right\vert ^{2}.
\end{equation*}%
In the wavelet notation this is%
\begin{equation*}
A\int_{0}^{\infty }\int_{-\infty }^{\infty }\left\vert W_{S_{n}^{\alpha }(%
\frac{t}{2})}f(-x,s)\right\vert ^{2}s^{-2}dxds\leq \sum_{j,m}\left\vert
W_{S_{n}^{\alpha }(\frac{t}{2})}f(a^{j}bm,a^{j})\right\vert ^{2}.
\end{equation*}%
Taking into account that $S_{n}^{\alpha }(\frac{t}{2})$ is admissible, we
can apply (\ref{isometry}) and (\ref{sampl})\ to write the above inequality
as 
\begin{equation}
A\left\Vert f\right\Vert ^{2}\leq \sum_{j,m}\left\vert \left\langle
f,T_{a_{m,j}}D_{b_{m,j}}S_{n}^{\alpha }(\frac{t}{2})\right\rangle
\right\vert ^{2}.  \label{Aframe}
\end{equation}%
This is the left frame inequality. The sequence $%
\{T_{a_{m,j}}D_{b_{m,j}}S_{n}^{\alpha }(\frac{t}{2})\}_{j,m}$ is thus a
frame. $\Box $

In view of remark 5.1, the lattice result is a special case of our main
theorem.

\begin{corollary}
If $b\log a<\frac{4\pi }{2n+\alpha }$, then $\{T_{a^{j}bm}D_{a^{j}}S_{n}^{%
\alpha }(\frac{t}{2})\}_{j,m}$ is a frame of $H^{2}(\mathbf{U})$.
\end{corollary}

\begin{remark}
Theorem 5.1 parallels theorem 3.1 in \cite{GrLy}, which states that, in the
lattice case of the time-frequency plane, if the density of the lattice $%
\Lambda $ is $>n+1$ (or if the size of $\Lambda $ is $<(n+1)^{-1}$) then the
Gabor system $\{e^{2\pi i\lambda _{2}t}H_{n}(t-\lambda _{1}):\lambda
=(\lambda _{1},\lambda _{2})\in \Lambda \}$, where $H_{n}$ stands for the
Hermite function of order $n$, is a frame for $L^{2}(\mathbf{R})$. In
particular, if one is dealing with the Von Neumann lattice with parameters $%
a $ and $b$, the condition is $ab<(n+1)^{-1}$. Corollary 5.2 makes this
analogy even more explicit.
\end{remark}

\begin{remark}
It is interesting to notice that, since $Ber^{\alpha }$ $f(z)=\int_{0}^{%
\infty }t^{\alpha }e^{izt}(\mathcal{F}f)(t)dt,$ we clearly have 
\begin{equation*}
i^{k}Ber^{k+\frac{\alpha }{2}}f(z)=\left( \frac{d}{dz}\right) ^{k}Ber^{\frac{%
\alpha }{2}}f(z).
\end{equation*}
This allows to rewrite (\ref{Formula}) as 
\begin{equation*}
W_{S_{n}^{\alpha }(\frac{t}{2})}f(-x,s)=C\sum_{k=0}^{n}\frac{a_{k}}{i^{k}}%
s^{k+\frac{\alpha }{2}+\frac{1}{2}}\left( \frac{d}{dz}\right) ^{k}Ber^{\frac{%
\alpha }{2}}f(z)
\end{equation*}%
which is, in some sense, reminiscent of Proposition 3.2 in \cite{GrLy}.
\end{remark}

\begin{remark}
Observe that combining the Paley-Wiener with the Plancherel theorem, we have 
$\left\Vert f\right\Vert _{H^{2}(\mathbf{U})}=\left\Vert f\right\Vert
_{L^{2}(0,\infty )}$ and the results of Theorem 1 and 2 say also that we
have a frame for $L^{2}(0,\infty )$.
\end{remark}

\section{Further properties of the functions $S_{n}^{\protect\alpha }(t)$}

The notation 
\begin{equation*}
F(a,b;c;x)=\sum_{k=0}^{\infty }\frac{(a)_{k}(b)_{k}}{k!(c)_{k}}x^{k}\text{.}
\end{equation*}
for the hypergeometric function is used in this section. Rewriting $%
S_{n}^{\alpha }$ in this notation gives 
\begin{equation*}
S_{n}^{\alpha }(t)=C^{\alpha ,n}\left( \frac{1}{2}-it\right) ^{-\frac{\alpha 
}{2}-1}\text{ }F(-n,\frac{\alpha }{2}+1;\alpha +1;\frac{1}{\frac{1}{2}-it})
\end{equation*}
(observe that the infinite sum becomes a polynomial of order $n$, since $%
(-n)_{k}=0$ if $k>0$). Composing the functions $S_{n}^{\alpha }(t)$ with the
fractional linear transformation 
\begin{equation}
z=\frac{2t-i}{2t+i}  \label{flt}
\end{equation}
the result is 
\begin{equation}
S_{n}^{\alpha }(t)=\Gamma \left( \frac{\alpha }{2}+1\right) (1-z)^{\frac{%
\alpha }{2}+1}g_{n}^{\alpha }(z)  \label{s}
\end{equation}
where 
\begin{equation*}
g_{n}^{\alpha }(z)=\frac{(\alpha /2)_{n}}{n!}F(-n,\frac{\alpha }{2}%
+1;-n-\alpha /2+1;z)
\end{equation*}
is a polynomial in $z$ of degree $n$. This was pointed out in \cite{Shen}.
It was also shown that these polynomials satisfy the orthogonality 
\begin{equation}
\int_{\left| z\right| =1}g_{n}^{\alpha }(z)\overline{g_{n}^{\alpha }(z)}%
\left| 1-z\right| ^{\alpha }\frac{dz}{z}=0\text{ if }m\neq n  \label{ort}
\end{equation}
and therefore are orthogonal on the unit circle with respect to the weight 
\begin{equation}
w(z)=\sin ^{\alpha }\frac{\theta }{2}d\theta .  \label{measure}
\end{equation}
This fact implies many properties, since there exists a very rich theory for
orthogonal polynomials on the unit circle (see \cite{Simon} and references
therein and also chapter 8 of \cite{Ism}). For example, the general theory
assures that all the zeros of $g_{n}^{\alpha }(z)$ lay within the unit disc.

\begin{remark}
Setting $a=\frac{\alpha }{w}$ in Example 8.2.5 at \cite{Ism}, and using the
identity $1-e^{i\theta }=4\sin ^{2}\frac{\theta }{2}$ to write the measure
(8.2.21) as (\ref{measure}) we recognize that the polynomials $g_{n}^{\alpha
}(z)$ are, up to a normalization, a family of orthogonal polynomials on the
unit circle known as the \emph{circular Jacobi orthogonal polynomials.}
\end{remark}

\begin{remark}
From (\ref{ort})\ the functions $z^{-\frac{1}{2}}(1-z)^{\frac{\alpha }{2}%
}g_{n}^{\alpha }(z)$ are orthogonal on the circle and form a basis of the
Hardy space on the unit disc 
\begin{equation*}
H^{2}(\mathbf{D})=\{f:f\text{ is analytic in }\mathbf{D}\text{ and }%
\sup_{r<1}\int_{0}^{2\pi }\left\vert f(re^{it})\right\vert ^{2}dt<\infty \}
\end{equation*}%
it is also clear that $S_{n}^{\alpha }(t)$ are orthogonal on the real line
(the boundary of the upper half place). Since $g_{n}^{\alpha }(z)$ are the
circular Jacobi orthogonal polynomials, the basis functions $z^{-\frac{1}{2}%
}(1-z)^{\frac{\alpha }{2}}g_{n}^{\alpha }(z)$ are the circular Jacobi
orthogonal functions and it is therefore natural to call $S_{n}^{\alpha }(t)$
the \emph{rational Jacobi orthogonal functions.}
\end{remark}

\begin{remark}
From the general theory \cite[(8.2.10)]{Ism} follows that, if $\kappa _{n}$
is the leading coefficient of the polynomial, then the sequence of
polynomials $\{g_{n}^{\alpha }(z)\}$ satisfies a three term recurrence
relation 
\begin{equation*}
\kappa _{n}g_{n}^{\alpha }(0)g_{n+1}^{\alpha }(z)+\kappa
_{n-1}g_{n+1}^{\alpha }(0)zg_{n-1}^{\alpha }(z)=[\kappa _{n}g_{n+1}^{\alpha
}(0)+\kappa _{n+1}g_{n}^{\alpha }(0)z]g_{n}^{\alpha }(z)
\end{equation*}%
where $\kappa _{n}$ is the leading coefficient of the polynomial. From the
explicit representation of the polynomials $g_{n}^{\alpha }(z)$ it is easily
seen that 
\begin{equation*}
\kappa _{n}=\frac{(\frac{\alpha }{2}+1)_{n}}{n!};\text{ \ \ }\phi _{n}(0)=%
\frac{\alpha (\frac{\alpha }{2}+1)_{n-1}}{2n!}
\end{equation*}%
This three term recurrence relation provides a effective method for
computational purposes: To evaluate the functions $S_{n}^{\alpha }(t)$ it is
sufficient to combine this recurrence relation with formulas (\ref{flt}) and
(\ref{s}).
\end{remark}

\end{document}